\def\VersionDateTime{23/September/2015, 21:00 GMT+9:00. Version $1.2$}
\documentclass[12pt]{amsart}
\usepackage{amssymb,amsthm,amsmath,amscd,rsfs}
\usepackage{latexsym}

\newcommand{\QQ}{{\mathbb{Q}}}

\newcommand{\codim}{\operatorname{codim}}

\newcommand{\Spec}{\operatorname{Spec}}

\newcommand{\Proof}{{\sl Proof.}\quad}
\newcommand{\QED}{{\unskip\nobreak\hfil\penalty50\quad\null\nobreak\hfil
{$\Box$}\parfillskip0pt\finalhyphendemerits0\par\medskip}}

\newcommand{\Tr}{\operatorname{Tr}}

\newcommand{\pr}{\operatorname{pr}}

\title
[
Bogomolov conjecture
]{Geometric Bogomolov conjecture for
nowhere degenerate abelian varieties of dimension $5$
with trivial trace}
\author
{Kazuhiko Yamaki}
\date{\VersionDateTime}
\address
{Institute for Liberal Arts and Sciences,
Kyoto University, Kyoto, 606-8501, Japan}
\email{yamaki.kazuhiko.6r@kyoto-u.ac.jp}

\begin{document}

\theoremstyle{plain}
\newtheorem{Theorem}{Theorem}[section]
\newtheorem{Lemma}[Theorem]{Lemma}
\newtheorem{Proposition}[Theorem]{Proposition}
\newtheorem{Corollary}[Theorem]{Corollary}
\newtheorem{Main-Theorem}[Theorem]{Main Theorem}
\newtheorem{Theorem-Definition}[Theorem]{Theorem-Definition}
\theoremstyle{definition}
\newtheorem{Definition}[Theorem]{Definition}
\newtheorem{Remark}[Theorem]{Remark}
\newtheorem{Conjecture}[Theorem]{Conjecture}
\newtheorem{Claim}{Claim}
\newtheorem{Example}[Theorem]{Example}
\newtheorem{Key Fact}[Theorem]{Key Fact}
\newtheorem{ack}{Acknowledgments}       \renewcommand{\theack}{}
\newtheorem*{n-c}{Notation and convention}      
\newtheorem{citeTheorem}[Theorem]{Theorem}
\newtheorem{citeProposition}[Theorem]{Proposition}
\newtheorem{citeLemma}[Theorem]{Lemma}

\newtheorem{Step}{Step}

\renewcommand{\theTheorem}{\arabic{section}.\arabic{Theorem}}
\renewcommand{\theClaim}{\arabic{section}.\arabic{Theorem}.\arabic{Claim}}
\renewcommand{\theequation}{\arabic{section}.\arabic{Theorem}.\arabic{Claim}}

\def\Pf{\trivlist\item[\hskip\labelsep\textit{Proof.}]}
\def\endPf{\strut\hfill\framebox(6,6){}\endtrivlist}

\def\Pfo{\trivlist\item[\hskip\labelsep\textit{Proof of Proposition~\ref{ch-of-hyp}.}]}
\def\endPfo{\strut\hfill\framebox(6,6){}\endtrivlist}

\maketitle


\begin{abstract}
We prove
that the geometric Bogomolov conjecture
holds
for nowhere degenerate abelian varieties
of dimension $5$ with trivial trace.
By this result together with our previous work,
we see that the conjecture
holds for an abelian variety
such that the difference between the dimension
of its maximal nowhere degenerate abelian subvariety
and that of its trace equals $5$.
\end{abstract}

\section{Introduction}


\subsection{Geometric Bogomolov conjecture}

Let $k$ be an algebraically closed field
and let
$\mathfrak{B}$ be a normal projective variety.
Let $K$ be the function field of $\mathfrak{B}$
unless otherwise specified.
We fix an algebraic closure $\overline{K}$ of $K$.
Once we take an ample line bundle $\mathcal{H}$ on $\mathfrak{B}$,
we have the notion of heights on projective varieties over $\overline{K}$.
We refer to \cite{lang2} for details.

Let $A$ be an abelian variety over $\overline{K}$.
Let
$L$ be an even ample line bundle on $A$,
where
``even'' means that the pull-back 
of $L$
by the $(-1)$-times endomorphism is isomorphic to $L$.
Then we have the canonical height $\widehat{h}_{L}$ on $A$
associated to $L$,
which is known to be a semipositive quadratic form on $A \left( 
\overline{K} \right)$.
Let $X$ be a closed subvariety of $A$.
For any $\epsilon > 0$,
put
\[
X ( \epsilon ; L)
:=
\left\{
x \in X \left(
\overline{K} \right)
\left|
\ 
\widehat{h}_{L}
(x) \leq \epsilon
\right.
\right\}
.
\]
Then by \cite[Lemma~2.1]{yamaki5},
the property 
that $X ( \epsilon ; L)$ is dense in $X$ for any $\epsilon > 0$
does not depend on the choice of the even ample line bundle $L$.
We say that \emph{$X$ has dense small points}
if it has that property
(cf. \cite[Definition~2.2]{yamaki5}).

It would be natural to wish to characterize those closed subvarieties
which have dense small points.
In the arithmetic setting,
some characterization has been established.
Indeed,
when $K$ is a number field,
Zhang has proved
in \cite{zhang2} that a closed subvariety has dense small points if and
only if it is a torsion subvariety,
i.e.,
the translate
of an abelian subvariety by a torsion point.
In \cite{moriwaki5},
Moriwaki has generalized
Zhang's theorem
to the case where 
$K$ is a finitely generated field over $\QQ$,
with respect arithmetic heights
introduced by himself.
In the case of function fields 
with respect to classical heights, however,
the characterization problem is still open.
Unlike
the arithmetic cases,
some abelian varieties
over function fields
actually
have closed subvarieties which are not
torsion subvarieties but have dense small points.
Taking into account that fact,
we
have proposed
in \cite{yamaki5}
the following conjecture,
called the \emph{geometric Bogomolov conjecture}.
\begin{Conjecture} [Conjecture~2.9 of 
\cite{yamaki5}] \label{conj:GBCforAVint}
Let $A$ be an abelian variety over $\overline{K}$.
Let $X$ be a closed subvariety of $A$.
Then $X$ has dense small points if and only if 
$X$ is a special subvariety.
\end{Conjecture}

Here,
we recall the definition of special subvarieties
defined in \cite{yamaki5}.
To define them,
we use the $\overline{K}/k$-trace of $A$.
A pair $\left( \widetilde{A}^{\overline{K}/k}, \Tr_A \right)$
consisting
of an abelian variety $\widetilde{A}^{\overline{K}/k}$ over $k$ and a homomorphism
$\Tr_A : \widetilde{A}^{\overline{K}/k} \otimes_{k} 
\overline{K} \to A$
of abelian varieties over $\overline{K}$
is called a \emph{$\overline{K}/k$-trace} of $A$ if
it has the following property:
For any abelian variety $\widetilde{B}$ over $k$
and a homomorphism $\phi : \widetilde{B} \otimes_{k} 
\overline{K}
\to A$,
there exists a unique homomorphism $\Tr (\phi) :
\widetilde{B} \to 
\widetilde{A}^{\overline{K}/k} $ such that
$\phi$ factors as
$\phi = \Tr_A \circ \left( \Tr (\phi) \otimes_k \overline{K} \right)$.
It is unique by the universal property,
and it is also known to exist.
We call $\Tr_A$ the \emph{$\overline{K}/k$-trace homomorphism} of $A$.
We refer to \cite{lang1} for more details.
Then,
a closed subvariety $X$ of $A$ is said to be \emph{special}
if there exists a closed subvariety $\widetilde{Y}$ of 
$\widetilde{A}^{\overline{K}/k}$
and a torsion subvariety $T$ of $A$ such that 
$X = \Tr_A \left( \widetilde{Y} \otimes_k \overline{K} \right)
+ T$.

Remark that
one easily sees that 
the ``if'' part of the geometric Bogomolov conjecture holds
(cf. \cite[Corollary~2.8]{yamaki5}),
and
the essential part
is the ``only if'' part.

Although
Conjecture~\ref{conj:GBCforAVint}
is still open in full generality,
there are some partial answers.
In \cite{gubler2},
Gubler has proved the conjecture holds
under the assumption that $A$ is totally degenerate
at some place.
In \cite{yamaki6},
we have generalized Gubler's work
and
have
shown that the conjecture for $A$ is reduced to the conjecture for
its maximal nowhere degenerate abelian subvariety $\mathfrak{m}$,
where ``being nowhere degenerate'' means ``having
everywhere good reduction''
and ``maximal'' means ``maximal with respect to the inclusion''
(cf. \cite[Definition~7.10]{yamaki6}).
It is known that 
there exists such an $\mathfrak{m}$ uniquely.
In \cite{yamaki7},
we have proved the following theorem,
where
$\mathfrak{t}$ is the image of the $\overline{K}/k$-trace
homomorphism $\Tr_A$  of $A$,
and it
is known to be contained in
$\mathfrak{m}$
(cf. \cite[Lemma~7.8~(2)]{yamaki6}).

\begin{citeTheorem} [cf. Theorem~1.5 of \cite{yamaki7}]
\label{thm:yamaki72int}
With the notation above,
the following are equivalent to each other:
\begin{enumerate}
\renewcommand{\labelenumi}{(\alph{enumi})}
\item
The geometric Bogomolov conjecture holds for $A$;
\item
The geometric Bogomolov conjecture holds for $\mathfrak{m}/\mathfrak{t}$.
\end{enumerate}
\end{citeTheorem}

Remark that
in a recent paper \cite{yamaki8},
we have shown that if 
$X$ is a closed subvariety
having dense small points such that
$\dim (X) = 1$ or $\codim (X , A) = 1$,
then it is a special subvariety.
That result will be used in the arguments in this note.

\subsection{Results}
In this note, we prove that
Conjecture~\ref{conj:GBCforAVint}
holds for nowhere degenerate abelian varieties
of dimension $5$ with trivial $\overline{K}/k$-trace:

\begin{Theorem} \label{thm:main1int}
Let $A$ be a nowhere degenerate abelian variety over $\overline{K}$
of dimension $5$ with trivial $\overline{K}/k$-trace.
Let $X$ be a closed subvariety of $A$.
Suppose that $X$ has dense small points.
Then $X$ is a torsion subvariety.
\end{Theorem}

Theorem~\ref{thm:main1int}
contributes the geometric Bogomolov conjecture
not only for nowhere degenerate abelian varieties
but also for more general abelian varieties.
Let $A$ be any abelian variety over $\overline{K}$.
Let $\mathfrak{m}$ 
and $\mathfrak{t}$
be 
as above.
Since $\mathfrak{m}/\mathfrak{t}$
is a nowhere degenerate abelian variety with trivial $\overline{K}/k$-trace
(cf. \cite[Lemma~7.8~(2)]{yamaki6} and \cite[Remark~5.4]{yamaki7}),
the following theorem is deduced from
Theorem~\ref{thm:main1int}
and
Theorem~\ref{thm:yamaki72int}.

\begin{Theorem} \label{thm:main2int}
Let $A$, $\mathfrak{m}$ and $\mathfrak{t}$ be as above.
Assume that $\dim ( \mathfrak{m} / \mathfrak{t} ) = 5$.
Then the geometric Bogomolov conjecture holds for $A$.
\end{Theorem}

In \cite{yamaki8},
we show that the geometric Bogomolov conjecture holds for
$A$ 
with $\dim ( \mathfrak{m} / \mathfrak{t} ) \leq 3$.
It follows 
by Theorem~\ref{thm:main2int}
that
the conjecture holds for 
$A$ such that $\dim ( \mathfrak{m} / \mathfrak{t} ) \leq 5$
and $\dim ( \mathfrak{m} / \mathfrak{t} ) \neq 4$.
At this moment,
we do not have an idea for $\dim ( \mathfrak{m} / \mathfrak{t} ) = 4$.

The proof of Theorem~\ref{thm:main1int}
relies
heavily on the results in \cite{yamaki8},
that is,
the non-density of small points
on closed subvarieties of dimension $1$ or of codimension $1$.
The key tools to deduce the theorem from
this non-density result are the morphisms
$X^{m} \to A^{m-1}$ given by 
$(x_1 , \ldots , x_m) \mapsto (x_1 - x_2 , x_2 - x_3 , \ldots ,
x_{m-1} - x_m )$ for $m = 2 ,3$.

\subsection*{Acknowledgments}
This work was initiated during my visit
to the university of Regensburg in March--April 2015,
which was supported by the SFB Higher Invariants.
I thank Professor Walter Gubler for inviting me 
and for his hospitality.
I thank him also for his comments.
This work was partly supported by KAKENHI 26800012.

\section{Proof}

In this article,
a \emph{variety} means a non-empty
geometrically integral scheme which is separated 
and of finite type over a field.

\subsection{Difference morphism}
In this subsection, let $A$ be an abelian variety over an
algebraically closed field $\mathfrak{K}$.
Let $X$ be a closed subvariety of $A$.
For any $x \in X ( \mathfrak{K} )$,
let $G(x)$
be the closed subset of $A$ defined by
\[
G (x)
:=
\{
a \in A \mid a+x \in X
\}
.
\]
For $x_1 , \ldots , x_m \in X ( \mathfrak{K} )$,
we set $G (x_1 , \ldots , x_m)
:= G (x_1) \cap \cdots \cap G (x_m)$.
Let $G_X$ be the stabilizer of $X$.
Then
$G_X = \bigcap_{x
\in X ( \mathfrak{K} )}
G (x)$.

For any integer $m$ with $m \geq 2$,
let $\alpha_{X^m} :
X^{m}
\to 
A^{m-1}$
be the morphism given by
\begin{align*}
\alpha_{X^m}
:
(x_1 , \ldots , x_m)
\mapsto
(x_1 - x_2 , x_2 - x_3 , \ldots , x_{m-1} - x_m)
.
\end{align*}
These morphisms $\alpha_{X^m}$ for $m \geq 2$, called difference morphisms,
will play important roles.

The following assertion is stated in the proof of
\cite[Lemma~3.1]{zhang2}.
It can be also found at \cite[Claim~6.2.2.1]{KMY}
with a proof.

\begin{citeLemma} \label{lem:G-fiber}
For any $( x_1 , \ldots , x_m )
\in X^{m } ( \mathfrak{K} )$,
we have
\[
( \alpha_{X^{m}} )^{-1}
(
\alpha_{X^{m}}
(
x_1 , \ldots , x_m
)
)
=
\{ (x_1 + a , \ldots , x_m + a)
\mid
a \in
G (x_1 , \ldots , x_m)
\}
.
\]
In particular,
$\dim ( ( \alpha_{X^{m}} )^{-1}
(
\alpha_{X^{m}}
(
x_1 , \ldots , x_m
)
) ) = \dim (G (x_1 , \ldots , x_m))$.
\end{citeLemma}

We set
\[
d_{m} (X) :=
\min \{ \dim (G (x_1 , \ldots , x_m) )
\mid x_1 , \ldots , x_m \in X ( \mathfrak{K})
\}
,
\]
where 
``$\dim$'' means the maximum of the dimensions
of the irreducible components.
Remark that $d_{m} (X) \leq \dim ( G (x_1 , \ldots , x_m) )$
for any $(x_1 , \ldots , x_m) \in X^{m} ( \mathfrak{K}
)$.
Remark also
that
$d_{m} (X) = \dim ( G (x_1 , \ldots , x_m) )$
for general $(x_1 , \ldots , x_m) \in X^{m} ( \mathfrak{K}
)$
by Lemma~\ref{lem:G-fiber}.

The following lemma is the key in our proof of the theorem.

\begin{Lemma} \label{lem:dmdm+1}
Suppose that $d_{m} (X) > \dim (G_X)$.
Then $d_{m+ 1 } (X) < d_{m} (X)$.
\end{Lemma}

\Proof
Fix an $(x_1 , \ldots , x_m) \in X ( \mathfrak{K} )$ with
$\dim ( G (x_1 , \ldots , x_m)) = d_m (X)$.
Let $Z_1 , \ldots , Z_p$ be the irreducible components
of $G (x_1 , \ldots , x_m)$
of dimension $d_m (X)$.
For any $j = 1 , \ldots , p$, 
let $\psi_j : Z_j \times X \to A$ be the restriction of the addition
$A^{2} \to A$
and
set $W_j := \psi_j (Z_j \times X)$.
Since $\dim (Z_j) = d_m (X) > \dim (G_X)$,
we have
$Z_j \nsubseteq G_X$.
It follows that
$W_j \nsubseteq X$,
and hence
$\psi_j^{-1} ( X ) \subsetneq Z_j \times X$.
Therefore,
$\pr_2 \left( ( Z_j \times X ) 
\setminus \psi_j^{-1} ( X ) \right)$,
where $\pr_{2} : A \times X \to X$ is the second projection,
is a dense open subset of $X$.
Thus $U := \bigcap_{j = 1}^{p} \pr_X \left( ( Z_j \times X ) 
\setminus \psi_j^{-1} ( X ) \right)$
is 
a dense open subset of $X$.

Take any $x_{m+1} \in U ( \mathfrak{K} )$.
Then for any $j = 1 , \ldots , p$,
there exists $a_j \in Z_j ( \mathfrak{K}
)$ such that $(a_j , x_{m+1}) \notin \psi_j^{-1} ( X )$,
which
means
$a_j \notin G (x_{m+1})$.
Thus
for any $j = 1 , \ldots , p$,
we have $G ( x_{m+1} ) \nsupseteq Z_j$.
It follows that
$G (x_1 , \ldots , x_m , x_{m+1}) = G (x_1 , \ldots , x_m)
\cap G (x_{m+1})$ does not contain any irreducible component
of $G (x_1 , \ldots , x_m)$
of dimension $d_m (X)$,
which implies that $\dim ( G (x_1 , \ldots , x_m , x_{m+1}) )
< d_m (X)$.
This shows $d_{m+1} (X) < d_m (X)$.
\QED

The following lemma
characterizes translates of abelian subvarieties in terms of the
difference morphism.

\begin{Lemma} \label{lem:alpha2}
Suppose that
$\dim \left( \alpha_{X^2} (X^{2}) \right) = \dim (X) $.
Then $X$ is a translate of an abelian subvariety.
\end{Lemma}

\Proof
Fix a closed point $x_0 \in X$
and set $Z = X - x_0$.
It then suffices to show that $Z - Z \subset Z$.
First, we note
that
$Z - Z = X - X  = \alpha_{X^{2}} (X^{2}) $
and that
$Z = \alpha_{X^2} ( X \times \{ x_0 \} )
\subset \alpha_{X^2} (X^{2}) $.
Since $\dim (\alpha (X^{2}) ) = \dim (X) = \dim (Z)$
and $\alpha (X^{2})$ is irreducible,
it follows that $Z = \alpha (X^{2})$,
Thus $Z - Z = \alpha_{X^{2}} (X^{2}) 
= Z$
as required.
\QED

\subsection{Proof of the main theorem}

For the proof of Theorem~\ref{thm:main1int},
we show three lemmas.
First two ones would be
more or less known facts.

\begin{Lemma} \label{lem:finallyadded}
Let $A$ be an abelian variety over $\overline{K}$
and let $X$ be a closed subvariety of $A$.
Suppose that $X$ is a translate of an abelian subvariety
and has dense small points.
Then it is a torsion subvariety.
\end{Lemma}

\Proof
We take an abelian subvariety $G$ of $A$ and an
$a_0 \in A \left( \overline{K} \right)$
such that $X = G + a_0$.
Let $\phi : A \to A / G$ be the quotient homomorphism.
Then $\phi (X) = \{ \phi (a_0) \}$.
By \cite[Lemma~2.1]{yamaki5},
this has dense small points,
which means that $\phi (a_0)$ is a torsion point.
Therefore,
there exists a torsion point $\tau \in A \left( \overline{K} \right)$
such that $\phi ( \tau ) = \phi (a_0)$
(cf.  the proof of \cite[Lemma~2.10]{yamaki5}).
Thus $X = G + \tau$,
which completes the proof of the lemma.
\QED

\begin{Lemma} \label{lem:dim=1ASV}
Let $A$ be a nowhere degenerate abelian variety over $\overline{K}$
with trivial
$\overline{K}/k$-trace.
Then $A$ does not have an abelian subvariety of dimension $1$.
\end{Lemma}

\Proof
To show the lemma by contradiction,
suppose that
there exists an abelian subvariety $E$ of $A$ 
with $\dim (E) = 1$.
By \cite[Lemma~5.5]{yamaki8},
$E$ is nowhere degenerate and has trivial $\overline{K}/k$-trace.
By \cite[Proposition~2.4]{yamaki7},
there exist a finite extension $K'$ of $K$ in $\overline{K}$,
an open subset $\mathfrak{U}$ of $\mathfrak{B}'$
with $\codim \left( \mathfrak{B}' \setminus \mathfrak{U},
\mathfrak{B}' \right) \geq 2$
where $\mathfrak{B}'$ is the normalization of $\mathfrak{B}$
in $K'$,
and
an abelian scheme $f : \mathscr{E} \to \mathfrak{U}$
such that $\mathscr{E} \times_{\mathfrak{U}} \Spec \left( 
\overline{K}
\right) \cong E$.
Since the moduli space of elliptic curves is 
the affine line,
this $\mathscr{E} \to \mathfrak{U}$
gives rise to a morphism from $\mathfrak{U}$ to the
affine line.
Since $E$ has trivial $\overline{K}/k$-trace,
this morphism is non-constant.
Thus we obtain a non-constant regular function $f$ on $\mathfrak{U}$.
Since $\mathfrak{B}'$ 
is normal and $\codim \left( \mathfrak{B}' \setminus \mathfrak{U},
\mathfrak{B}' \right) \geq 2$,
$f$ extends to a regular function on $\mathfrak{B}'$.
However, that is impossible since $\mathfrak{B}'$ is projective.
Thus the lemma holds.
\QED

The other lemma is
a restatement of the non-density results
in \cite{yamaki8}
which is modified for our setting.

\begin{Lemma} \label{lem:dim-codim-neq1}
Let $A$ be a nowhere degenerate abelian variety over $\overline{K}$
with trivial
$\overline{K}/k$-trace
and
let $X$ be a closed subvariety of $A$.
Suppose that $X$ has dense small points.
Then $\dim (X) \neq 1$ and $\codim ( X , A) \neq 1$.
\end{Lemma}

\Proof
By
\cite[Proposition~5.6]{yamaki8},
$X$ is not a divisor, and thus
$\codim (X) \neq 1 $.
We show that $\dim (X) \neq 1$
by contradiction.
Suppose that $\dim (X) = 1$.
By \cite[Theorem~5.10]{yamaki8},
$X$ is a special subvariety of $A$.
Since $A$ has trivial $\overline{K}/k$-trace,
that implies that $X$ is a torsion subvariety of dimension $1$,
and thus
there exists an abelian subvariety $E$ of $A$ 
with $\dim (E) = 1$.
However, that contradicts Lemma~\ref{lem:dim=1ASV}.
\QED

We now restate the main theorem and prove it.

\begin{Theorem}[Theorem~\ref{thm:main1int}] \label{thm:main1}
Let $A$ be a nowhere degenerate abelian variety 
of dimension $5$ with trivial
$\overline{K}/k$-trace.
Let $X$ be a closed subvariety of $A$.
Suppose that $X$ has dense small points.
Then $X$ is a torsion subvariety.
\end{Theorem}

\Proof
The case of $\dim (X) = 5$ is trivial,
and the case of $\dim (X) = 0$ is classically known
(cf. \cite[Remark~7.4]{yamaki6}).
Further
by Lemma~\ref{lem:dim-codim-neq1},
we have $\dim (X) \neq 1$
and
$\dim (X) \neq 4$.
Thus we may assume that $\dim (X) = 2,3$.

To ease notation,
we write $\alpha_m$ for $\alpha_{X^{m}}$.
We set $Y_{2} := \alpha_{2} (X^{2})
\subset A$.
Since $X^{2}$ has dense small points,
so does $Y_2$.
By Lemma~\ref{lem:dim-codim-neq1},
we have $\dim (Y_2) \neq 4$.

First, we consider the case of $\dim (X) = 3$.
We prove $\dim (Y_2) \neq 5$ then.
To argue by contradiction,
suppose that $\dim (Y_2) = 5$.
Then $\alpha_2 : X^{2} \to A$ is surjective.
Since $\dim (X^{2} ) - \dim (A) = 1$,
we have $d_{2} (X) = 1$ (cf. Lemma~\ref{lem:G-fiber}).
Remark in particular that $\dim (G_X) \leq 1$.
By Lemma~\ref{lem:dim=1ASV},
we find
$\dim (G_X) = 0$.
By Lemma~\ref{lem:dmdm+1},
it follows that $d_{3} (X) = 0$.
This means that
$\alpha_3 : X^{3} \to A^{2}$ is generically finite,
and hence $\dim 
\left( \alpha_3 \left( X^{3} \right) \right) = 3 \dim (X) = 9$.
Thus
$\alpha_3 \left( X^{3} \right) $ is an
effective divisor on $A$.
On the other hand,
since $X$ has dense small points, so does $X^{3}$,
and hence $\alpha_3 \left( X^{3}  \right)$ has dense small points.
However,
since
$A^{3}$ is nowhere degenerate and has trivial $\overline{K}/k$-trace
(cf. \cite[Lemma~5.5]{yamaki8}), 
that contradicts 
Lemma~\ref{lem:dim-codim-neq1}.
Thus we obtain $\dim (Y_2) \neq 5$.

Since $\dim ( Y_2 ) \geq \dim (X) = 3$,
it follows that
$\dim (Y_2) = 3$.
Then
Lemma~\ref{lem:alpha2}
shows
that
$X$ is a translate of an abelian subvariety.
Since $X$ has dense small points,
it follows 
by Lemma~\ref{lem:finallyadded}
that $X$ is a torsion subvariety.
Thus we obtain the theorem in the case of $\dim (X) = 3$.

Next, we consider the case of $\dim (X) = 2$.
Since $\dim (Y_2) \neq 4$,
we then have
$2 \leq \dim (Y_2) \leq 3$.
In fact, we prove $\dim (Y_2) =2$.
To argue by contradiction,
suppose that $\dim (Y_2) = 3$.
Since $Y_2$ is a closed subvariety of $A$
having dense small points with $\dim (Y_2) = 3$,
it follows from what we have shown above that $Y_2$ is a torsion subvariety.
Since $0 \in Y_2$, $Y_2$ is an abelian subvariety.
By the Poincar\'e complete reducibility
theorem,
therefore,
there exists an abelian subvariety $A'$ of $A$ of dimension $2$ such that
$Y_2 + A' = A$.
Let $\phi : A \to A/A'$ be the quotient homomorphism.
Then $A/A'$ is an abelian variety of dimension $3$,
and furthermore
by \cite[Lemma~5.5]{yamaki8},
it is nowhere degenerate and has trivial $\overline{K}/k$-trace.
Since
$Y_2 \cap A'$ is finite
and since some translate of $X$ lies in $Y_2$,
we see that
$\dim (\phi (X)) = \dim (X) = 2$.
Thus
$\phi (X)$ is a divisor on $A/A'$.
On the other hand, since $X$ has dense small points,
so does $\phi (X)$.
However, that is a contradiction by
Lemma~\ref{lem:dim-codim-neq1}.

Thus we have $\dim (Y_2) = 2$.
By Lemma~\ref{lem:alpha2},
it follows that
$X$ is a translate of an abelian subvariety.
Since $X$ has dense small points,
Lemma~\ref{lem:finallyadded}
concludes that $X$ is a torsion subvariety.
Thus we complete the proof.
\QED

Finally,
as is noted in the introduction,
Theorem~\ref{thm:main2int}
is deduced
from Theorem~\ref{thm:main1}.


\renewcommand{\thesection}{Appendix} 
\renewcommand{\theTheorem}{A.\arabic{Theorem}}
\renewcommand{\theClaim}{A.\arabic{Theorem}.\arabic{Claim}}
\renewcommand{\theequation}{A.\arabic{Theorem}.\arabic{Claim}}
\renewcommand{\theProposition}
{A.\arabic{Theorem}.\arabic{Proposition}}
\renewcommand{\theLemma}{A.\arabic{Theorem}.\arabic{Lemma}}
\setcounter{section}{0}
\renewcommand{\thesubsection}{A.\arabic{subsection}}



\small{

}

\end{document}